\documentclass[twoside,english,american]{elsarticle}
\usepackage[T1]{fontenc}
\usepackage[utf8]{inputenc}
\usepackage{float}
\usepackage{amsmath}
\usepackage{amsthm}
\usepackage{amssymb}
\usepackage{graphicx}
\usepackage{xcolor}

\makeatletter
	
\providecommand{\tabularnewline}{\\}

\theoremstyle{plain}
\newtheorem{thm}{\protect\theoremname}
\theoremstyle{definition}
\newtheorem{defn}[thm]{\protect\definitionname}
\theoremstyle{definition}
\newtheorem{example}[thm]{\protect\examplename}
\ifx\proof\undefined
\newenvironment{proof}[1][\protect\proofname]{\par
	\normalfont\topsep6\p@\@plus6\p@\relax
	\trivlist
	\itemindent\parindent
	\item[\hskip\labelsep\scshape #1]\ignorespaces
}{%
	\endtrivlist\@endpefalse
}
\providecommand{\proofname}{Proof}
\fi

\journal{Int. J. Approx. Reason.}

\@ifundefined{showcaptionsetup}{}{%
 \PassOptionsToPackage{caption=false}{subfig}}
\usepackage{subfig}
\makeatother

\usepackage{babel}
\addto\captionsamerican{\renewcommand{\definitionname}{Definition}}
\addto\captionsamerican{\renewcommand{\examplename}{Example}}
\addto\captionsamerican{\renewcommand{\theoremname}{Theorem}}
\addto\captionsenglish{\renewcommand{\definitionname}{Definition}}
\addto\captionsenglish{\renewcommand{\examplename}{Example}}
\addto\captionsenglish{\renewcommand{\theoremname}{Theorem}}
\providecommand{\definitionname}{Definition}
\providecommand{\examplename}{Example}
\providecommand{\theoremname}{Theorem}

\usepackage{etoolbox}
\makeatletter
\patchcmd{\ps@pprintTitle}
  {Preprint submitted to}
  {Preprint accepted to publication in}
  {}{}
\makeatother

\begin{document}
\begin{frontmatter}{}

\title{On the Derivation of Weights from Incomplete Pairwise Comparisons
Matrices via Spanning Trees with Crisp and Fuzzy Confidence Levels\tnoteref{t1}}

\tnotetext[t1]{These authors contributed equally to this work.}

\author[rvt]{Jiri Mazurek\corref{cor2}}

\ead{mazurek@opf.slu.cz}

\author[focal]{Konrad~Ku\l akowski \corref{cor1}}

\ead{konrad.kulakowski@agh.edu.pl}

\cortext[cor1]{Corresponding author}

\address[rvt]{Department of Informatics and Mathematics, Silesian University in
Opava, School of Business Administration in Karvina, Univerzitni namesti
1934/3, Czech Republic}

\address[focal]{Department of Applied Computer Science, AGH University of Science
and Technology, al. Mickiewicza 30, 30-159 Kraków, Poland}
\begin{abstract}
In this paper, we propose a new method for the derivation of a priority
vector from an incomplete pairwise comparisons (PC) matrix. We assume
that each entry of a PC matrix provided by an expert is also evaluated
in terms of the expert's confidence in a particular judgment. Then,
from corresponding graph representations of a given PC matrix, all
spanning trees are found. For each spanning tree, a unique priority
vector is obtained with the weight corresponding to the confidence
levels of entries that constitute this tree. At the end, the final
priority vector is obtained through an aggregation of priority vectors
achieved from all spanning trees. Confidence levels are modeled by
real (crisp) numbers and triangular fuzzy numbers. Numerical examples
and comparisons with other methods are also provided. Last, but not
least, we introduce a new formula for an upper bound of the number
of spanning trees, so that a decision maker gains knowledge (in advance)
on how computationally demanding the proposed method is for a given
PC matrix. 
\end{abstract}
\begin{keyword}
pairwise comparisons \sep fuzzy numbers \sep priority vector \sep
spanning tree \sep multiple-criteria decision making 
\end{keyword}

\end{frontmatter}{}

\section{Introduction}

Pairwise comparisons (PCs) constitute a fundamental part of many multiple-criteria
decision making (MCDM) methods, such as the AHP/ANP, BWM, ELECTRE,
MACBETH, PAPRIKA and PROMETHEE etc. \citep{Kulakowski2020utahp,Brunelli2019ambw,Greco2016mcda,Mardani2015mcdm,Saaty2004dmta,Vaidya2006ahpa,Zavadskas2014soas}.

Usually, it is assumed that all pairwise comparisons by experts are
performed and available to a decision maker. However, due to a lack
of knowledge, time pressure, uncertainty or other factors, pairwise
comparisons might not be complete. The problem of the derivation of
a priority vector from incomplete PC matrices has been studied for
almost four decades. Probably the first study on incomplete PC matrices
came in 1987 from Weiss and Rao \citep{Weiss1987adfl} in the context
of large-scale systems. Soon afterwards, the 'standard' method for
obtaining a priority vector from an incomplete PC matrix was proposed
by Harker in \citep{Harker1987amoq}. Harker's completion method is
based upon the graph-theoretic structure of the PC matrix and the
gradient of the right Perron vector. It can be used with the geometric
mean leading to the geometric mean method (GMM) for incomplete PC
matrices proposed by Ku\l akowski \citep{kulakowski2020otgm}. This
approach is equivalent to the logarithmic least-square method (LLSM)
for incomplete PC matrices proposed and developed earlier by Tone and
Bozóki et al. \citep{Tone1993llsm,Bozoki2010ooco}. Ram{\'\i}k
\citep{Ramik2020pcmt} proposed estimation of missing matrix elements
via spanning trees of an induced graph of an extended matrix while
maximally preserving consistency. The method of determining missing
values based on inconsistency was also proposed by Alonso et al. \citep{Alonso2008acbp}.
Several recent studies were aimed at theoretical considerations regarding
the geometric mean method. Lundy et al. \citep{Lundy2016tmeo} showing
that the ranking (a priority vector) based on spanning trees of a
PC matrix is equivalent to the ranking obtained by the geometric mean
method. Later on, Bozóki and Tsyganok \citep{Bozoki2019tlls} proved
that the spanning trees method and the logarithmic least squares method
for incomplete PC matrices yield the same result as well.

By far the most popular methods for the derivation of a priority vector
are the aforementioned eigenvalue method and geometric mean method
(the least squares method). In 2012, Siraj et al. \citep{Siraj2012east}
proposed the spanning trees method, originally for complete PC matrices.
The method enumerates all spanning trees of a graph representation
of a given PC matrix to find priority vectors corresponding to each
spanning tree, and, after aggregation, the final priority vector is
achieved. The virtue of this approach lies in the fact that, from
each spanning tree, a unique priority vector is obtained and that
pairwise comparisons included in one spanning tree are necessarily
consistent. However, the obvious drawback of this method is the large
number (given by Cayley's theorem) of spanning trees when the order
of a PC matrix grows. Nevertheless, for incomplete matrices the number
of spanning trees might be significantly reduced and the method can
be easily applicable.

Another problem in pairwise comparisons relates to the problem of
uncertainty associated with experts' judgments. Uncertainty in the
framework of PCs has been studied by many authors who proposed various
sophisticated methods and approaches based on interval numbers, fuzzy
sets, intuitionistic fuzzy sets, fuzzy hesitant sets, rough sets,
linguistic variables, Z-numbers, and so on \citep{Azadeh2013zahp,WilliamWest2021dtfw,Janicki2016oawr,Benitez2019moup,Dede2021afoi,Deng1999mawf,Durbach2012muim,Ramik2020pcmt,Yang2013hthc,Zadeh2011anoz},
to name a few, however, a comprehensive review on uncertainty in the
pairwise comparisons method is beyond the scope of this paper.

Therefore, this paper aims to propose a novel procedure for the derivation of a priority vector (a vector of weights of compared objects) from an incomplete PC matrix which is based on the graph representation
of a PC matrix and its set of spanning trees, which allows modeling
of uncertainty using confidence levels. The method is based
on elaboration and synthesis of the aforementioned studies, where
the application of spanning trees was inspired by \citep{Siraj2012east}
and the use of confidence levels by \citep{Zadeh2011anoz}. More specifically,
we introduce two similar mathematical models with different levels
of uncertainty. In both models, confidence levels for each pairwise
comparison provided by an expert are provided, creating a confidence
matrix. In the first model, the confidence levels are expressed (coded)
by crisp numbers corresponding to a given linguistic scale. For instance,
an expert can have ``strong confidence'' or ``absolute confidence''
in his/her judgments. Hence, the input of the first model consists
of a real reciprocal PC matrix and a corresponding real symmetric
confidence matrix, and the output is a real priority vector. In our
second model, we incorporate triangular fuzzy numbers both for the
preferences and confidences. Preferences form a reciprocal fuzzy PC
matrix and the confidence values are arranged into a symmetric fuzzy
confidence matrix. Both matrices constitute the input of the model,
while the output is a fuzzy priority vector. In both models, confidence
values are applied to assess the reliability of a priority vector
inferred from each spanning tree, where reliability plays the role
of weights in the aggregation process by which a final priority vector
is obtained. In the presented model, we use triangular fuzzy numbers
due to their easy and intuitive interpretation, simple arithmetics,
and a large number of their previous applications in the literature.
However, our approach can be easily extended to trapezoidal or other
types of fuzzy numbers as well.

It has several advantages. Firstly, it does not require a decision
maker to complete a PC matrix, which is always associated with a distortion
of original preferences. Secondly, it enables an expert to assign
confidence levels to his/her judgments, so the more confident preferences
(in a particular spanning tree) translate into the higher weight of
a corresponding priority vector in the final aggregation. Additionally,
incorporation of fuzzy numbers allows an expert to add another level
of uncertainty of his/her judgments dealing with pairwise comparisons.
Since the number of all pairwise comparisons grows quadratically with
the number of compared objects $n$, PC methods such as AHP are often
limited to small-scale problems, typically with $n$ not greater than
10. However, our approach based on incomplete pairwise comparisons
enables to solve problems with relatively high numbers of compared
objects as long as the amount of spanning trees induced by an incomplete
PC matrix remains manageable. That's why we provide an estimate for
the upper bound of the number of spanning trees as well. Finally,
the method is also applicable for complete PC matrices, however, due
to the fact that the number of spanning trees grows very quickly with
the size of a PC matrix, we recommend its use predominantly in the
former case.

The paper is organized as follows: Section 2 provides preliminaries
on applied concepts and notation, in Sections 3-4 methods for the
derivation of priority vectors are proposed with numerical examples,
and the Discussion with Conclusions close the article.

\section{Preliminaries}

\subsection{Multiplicative pairwise comparisons}

A pairwise comparison is a binary relation over a finite and discrete
set of objects $V=\lbrace V_{1},...,V_{n}\rbrace$ (usually formed
by alternatives, criteria, sub-criteria, etc.) such that $V\times V\rightarrow\mathbb{R}_{+}\cup\{0\}$.
The value 0 means the comparison is undefined. Further on, every pairwise
comparison $a_{ij}\in\mathbb{R}_{+}$ describes the degree of preference
(importance, etc.) of an object $V_{i}$ over an object $V_{j}$.

Pairwise comparisons are usually arranged into a square and positive
$n\times n$ matrix $A=(a_{ij})$ called a pairwise comparisons matrix.
If matrix contains zeros it is called incomplete. In our paper, we
limit ourselves to the multiplicative framework, however, multiplicative
pairwise comparisons can be easily transformed into additive or fuzzy
systems, and vice versa, see e.g. \citep{Cavallo2010cocp,Cavallo2012dwfa,Cavallo2018aguf}
or \citep{Ramik2020pcmt}. 
\begin{defn}
\label{def:A-matrix-recip} A matrix $A$ is said to be (multiplicatively)
reciprocal if: 
\begin{equation}
\forall i,j\in\{1,\ldots,n\}:a_{ij}\neq0\Rightarrow a_{ij}=\frac{1}{a_{ji}}
\end{equation}

and $A$ is said to be (multiplicatively) consistent if: 
\begin{equation}
\forall i,j,k\in\{1,\ldots,n\}:a_{ij},a_{jk},a_{ki}\neq0\Rightarrow a_{ij}\cdot a_{jk}\cdot a_{ki}=1.
\end{equation}
\end{defn}
The result of the pairwise comparisons method is a priority vector
(vector of weights) $w=(w_{1},...,w_{n})$ that assigns positive values
$w_{i}$ to each of the $n$ compared objects. According to the EVM
(the eigenvalue method) proposed by Saaty for complete PC matrices,
see \citep{Saaty1977asmf,Saaty2004dmta}, the vector $w$ is determined
as the rescaled principal eigenvector of a PC matrix $A$. Thus, assuming
that $Aw=\lambda_{\textit{max}}w$ the priority vector $w$ is given
as: 
\[
w=\left[\gamma w_{1},\ldots,\gamma w_{\textit{n}}\right]^{T},
\]
where $\gamma$ is a scaling factor. Usually, it is assumed that $\gamma=[\sum_{i=1}^{n}w_{i}]^{-1}$.
Extension of EVM for incomplete PC matrices has been defined by Harker
\citep{Harker1987amoq}. 

In the geometric mean method (GMM) introduced by Crawford and Williams
\citep{Crawford1985anot}, the weight of the i-th alternative is estimated
as the geometric mean of the i-th row of $A$. Thus, the priority
vector is given as: 
\begin{equation}
w=\left[\gamma\left(\prod_{r=1}^{n}a_{1r}\right)^{\frac{1}{n}},\ldots,\gamma\left(\prod_{r=1}^{n}a_{nr}\right)^{\frac{1}{n}}\right]^{T},\label{eq:GMM_def}
\end{equation}
where $\gamma$ is a scaling factor again. In the case of incomplete
matrices, one can use the extension of GMM defined by Ku\l{}akowski
\citep{kulakowski2020otgm} or equivalent Logarithmic Least-squares
Method (LLSM) \citep{Tone1993llsm,Bozoki2010ooco}.

Other, less frequently applied prioritization methods include direct
least squares, weighted least squares, row and/or column sums, etc.,
see e.g. \citep{Kulakowski2020utahp}.

Recently, an interesting method was proposed by Siraj et al. \citep{Siraj2012east,Kulakowski2020utahp}.
This method is based on spanning trees of a graph representation of
a given PC matrix, and the final priority vector is the arithmetic
mean of all spanning trees' priority vectors:

\begin{equation}
w=\frac{1}{\eta}\sum_{s=1}^{\eta}w(\tau_{s}),
\end{equation}

where $\eta=n^{n-2}$ denotes the number of spanning trees for a complete
PC matrix, see e.g. \citep{Cayley1889atot}, and $w(\tau_{s})$ denotes
priority vectors of individual spanning trees.

\subsection{Graph representation of a PC matrix and spanning trees}

Judgments contained in a PC matrix can be represented in a form of
a a graph in which vertices correspond to the alternatives while edges
denote individual pairwise comparisons, see e.g. \citep{Kulakowski2020utahp,Diestel2005gt,Siraj2012east,Kulakowski2018iito}. 
\begin{defn}
The undirected graph $P_{A}=(V,E)$ is said to be a graph of the PC
matrix $A=(a_{ij})$ if $V=\{V_{1},\ldots,V_{n}\}$ is the set of
vertices and $E\subseteq\{\{V_{i},V_{j}\}$ such that $\,\,V_{i},V_{j}\in V,\,\,i\neq j$
and $a_{ij}$ exists\footnote{For incomplete PC matrices, some comparisons may not be defined \citep{Kulakowski2020utahp}.}$\}$
is the set of edges. 
\end{defn}
Similarly, we can formally define the concept of a path between two
vertices. 
\begin{defn}
An ordered sequence of distinct vertices $p=V_{i_{1}},V_{i_{2}},\ldots,V_{i_{m}}$
such that $\left\{ V_{i_{1}},V_{i_{2}},\ldots,V_{i_{m}}\right\} \subseteq V$
is said to be a path between $V_{i_{1}}$ and $V_{i_{m}}$ with the
length $m-1$ in $P_{A}=(A,E)$ if $\{V_{i_{1}},V_{i_{2}}\}$, $\{V_{i_{2}},V_{i_{3}}\},\ldots$,$\{V_{i_{m-1}},V_{i_{m}}\}\in E$. 
\end{defn}
In a connected graph, there will be at least one path for each pair
of vertices. 
{
\begin{defn}
The irreducible matrix is one that cannot be transformed by the permutation
of rows and columns to form: 
\[
\left(\begin{array}{cc}
Q_{1} & 0\\
Q_{2} & Q_{3}
\end{array}\right),
\]
where $Q_{1},Q_{3}$ are square matrices, and $0$ denotes the entries
filled with zeros. 
\end{defn}
In the context of a PC matrix, zero means that the given comparison
is undefined, and a matrix having such comparisons will be called
incomplete. The connectivity of $P_{A}$ is equivalent to the irreducibility
of the matrix $A$ \citep[p. 185,  186]{Quarteroni2000nm}. It can
be shown that, in the case of incomplete PC matrices, the ranking
(a priority vector, vector of weights) can only be calculated if $P_{A}$
is connected, i.e. $A$ is irreducible\footnote{ Note that when the graph $P_{A}$ is undirected its matrix $A$ must
be symmetric. Thus, the reducibility means that there exists permutation
of rows and columns such that $A$ takes the form $A=\left(\begin{array}{cc}
Q_{1} & 0\\
0 & Q_{3}
\end{array}\right)$, where $Q_{1}$ and $Q_{3}$ are the square sub-matrices with the
dimensions $s\times s$ and $t\times t$ respectively. Assuming that
A corresponds to the comparisons of $s+t$ alternatives $V_{1},\ldots,V_{s+t}$
it is clear that for every $i,j$ such that $1\leq i\leq s<j\leq s+t$
the direct comparison $a_{ij}$ is $0$. This implies that non of
$V_{1},\ldots,V_{s}$ is compared to any of $V_{s+1},\ldots,V_{s+t}$.
Hence, based on $A$, no common ranking for $V_{1},\ldots,V_{s+t}$
can be calculated.} .
} 
\begin{defn}
The spanning tree $ST=(A,E')$ of $P_{A}$ is any connected subgraph
of $P_{A}$ such that $E'\subset E$ where the number of edges $\vert E'\vert=n-1$. 
\end{defn}
In the case of a complete PC matrix (the case of a complete graph),
the number of all spanning trees is given by the well-known Cayley
formula $n^{n-2}$ \citep{Cayley1889atot}, where $n$ is the order
of the PC matrix. However, this number is significantly reduced for
incomplete matrices (see e.g. Examples 1 and 3).
{
\begin{defn}
The Laplacian matrix $L(G)=(l_{ij})$ of a graph $G=(V,E)$ is a matrix
such that 
\[
l_{ij}=\begin{cases}
s_{i} & \text{if}\,i=j\\
0 & \text{if}\,i\neq j\,\,\text{and}\,\,\{v_{i},v_{j}\}\in E\\
-1 & \text{\text{if}\,i\ensuremath{\neq}j\,\,\text{and}\,\,\{\ensuremath{v_{i}},\ensuremath{v_{j}}\}\ensuremath{\notin}E}
\end{cases},
\]

where $V$ is a set of vertices $V=\{v_{1},\ldots,v_{n}\}$, $E$
is a set of edges in the form of pairs $\{v_{i},v_{j}\}$, where $v_{i},v_{j}\in V$,
and $s_{i}$ is the number of edges adjacent to $v_{i}$ i.e. $\{v_{i},v_{j}\}\in E$. 
\end{defn}
} 
The following theorems answer the question about the number of spanning
trees for incomplete graphs. 
\begin{thm}
\emph{Kirchhoff's matrix tree theorem} \citep{Kirchhoff1847udad}.
Let $G$ be a graph, and let $L(G)$ be the Laplacian matrix of $G$
and let $t(G)$ denote the number of spanning trees of a graph $G$.
Then $t(G)$ is equal to all cofactors of $L(G)$. 
\end{thm}
\begin{thm}
\emph{Kelmans and Chelnokov} \citep{Kelmans1974acpo}. Let $G$ be
a graph of $n$ points, let $L(G)$ be the Laplacian matrix of $G$
and let $0=\lambda_{1}\leq\lambda_{2}\leq...\leq\lambda_{n}$ denote
the eigenvalues of $L(G)$. Then: $t(G)=\frac{1}{n}\prod_{k=2}^{n}\lambda_{k}$. 
\end{thm}
In particular, if $G$ is a tree itself, then the number of spanning
trees $t(G)=1$; when $G$ is a cycle graph of $n$ points, then $t(G)=n$.

A spanning tree corresponding to a $n\times n$ PC matrix $A=(a_{ij})$
induces a unique ranking of compared objects since $(n-1)$ edges
of a spanning tree corresponds to $(n-1)$ equations for $n$ weights
of $n$ objects and along with the normality condition \textendash{}
the sum of weights being equal to 1 \textendash{} a system of $n$
equations is formed. Hence, a unique (normalized) priority vector
$w$ is assigned to each spanning tree.

\subsection{Pairwise comparisons' confidence and priority vectors' reliability}

According to \citep{Peterson1988cuat} or \citep{Wesson2009veoc},
the word 'confidence' is used to describe a person's strength of belief
about the accuracy or quality of a prediction, judgment, or choice,
and confidence can thus be described on a continuum ranging from total
certainty to complete doubt, though usually confidence is expressed
verbally on a discrete scale. 

For each $n\times n$ PC matrix $A=(a_{ij})$, let us define a
corresponding $n\times n$ ''confidence'' matrix $C=(c_{ij})$,
where $c_{ij}$ might be real (crisp) or fuzzy numbers, or possibly
other numbers, such as Z-numbers, such that $c_{ij}$ expresses the
confidence of an expert of a value of the pairwise comparison $a_{ij}$,
and it is assumed that $c_{ij}=c_{ji},\forall i,j$.

Further on, we assume that an expert can express his/her confidence
in PC judgments on a) a finite, discrete and strictly increasing cardinal
scale $S_{d}=\lbrace s_{1}^{(d)},...,s_{k}^{(d)}\rbrace$, where $s_{i+1}^{(d)}>s_{i}^{(d)},\forall i$,
b) a finite ordinal (linguistic) scale $S_{o}=\lbrace s_{1}^{(o)},...,s_{k}^{(o)}\rbrace$,
where $s_{i+1}^{(o)}\succ s_{i}^{(o)},\forall i$, or c) on a real
interval scale $S_{i}=[0,1]$.

Then we define the \textit{absolute reliability} $r$ of a priority
vector $w$ derived from a given spanning tree as follows: 
\begin{defn}
Let $A=(a_{ij})$ be a preference PC matrix, let $C=(c_{ij})$ be
a corresponding confidence matrix and let $E^{(ST)}=\lbrace[ij]\rbrace$
be a set of all edges of a spanning tree $ST$. Then the \textit{absolute
reliability} $r$ of a priority vector $w$ obtained from a spanning
tree $ST$ is given as follows: 
\begin{equation}
r(ST)=\lbrace\prod c_{ij}\vert[ij]\in E^{(ST)}\rbrace
\end{equation}
\end{defn}
Simply put, the value of the absolute reliability $r$ of a priority
vector obtained from a given spanning tree is a product of all $c_{ij}$
corresponding to the PCs' elements $a_{ij}$ (graph edges) of that
spanning tree. Next, we define the relative reliability $R$ of a
priority vector obtained from a given spanning tree. 

\begin{defn}
Let $r(ST)$ be the absolute reliability of a priority vector $w$
obtained from a spanning tree $ST$. Then the \textit{relative reliability}
$R$ of the priority vector $w$ obtained from a spanning tree $ST$ is given as follows: 
\begin{equation}
R(ST)=r(ST)/\sum_{i=1}^{N}r(ST_{i}),
\end{equation}
where $N$ is the number of all spanning trees corresponding to the
graph representation of a given PC matrix. 
\end{defn}

It should be noted that other aggregation functions, see Grabisch \citep{Grabisch2009eoma}, such as the minimum, can be applied for the calculation of the absolute reliability (5), though we consider the product to be the most natural. \footnote{Let's consider confidence values from Table 3 and a spanning tree with only two edges/preferences, with corresponding values of confidence 1 (low confidence) and 4 (absolute confidence). Then the absolute reliability of a priority vector acquired from this spanning tree is $ r = 1\cdot 4 = 4 $. Now consider another spanning tree with corresponding values of confidence 1 and 1. Then the absolute reliability $ r = 1\cdot 1 = 1 $. If we used, for example, the minimum function, both absolute reliabilites would be 1, neglecting the fact that an expert was more confident in the case of preferences for the first spanning tree, hence the corresponding priority vector shoud obtain the higher weight. }

As for the most suitable scale for the confidence values $c_{ij}$,
we suggest using the 5-item integer scale from 0 (zero confidence
in the judgment of $a_{ij}$) to 4 (absolute confidence about $a_{ij}$),
see Table \ref{tab:the-fuzzy-linguistic}. Similar scales have been
applied, for instance, in \citep{Lewin2015uqei}, or \citep{Li2020crid}.


\subsection{Fuzzy sets and triangular fuzzy numbers}

Fuzzy sets introduced by L. A. Zadeh \citep{Zadeh1965fs} provide
a convenient theoretical framework for uncertainty modeling. Here,
we provide its basic properties. 
\begin{defn}
A \textit{fuzzy set} $\tilde{U}$ is a pair $(U,\mu_{U}(x))$ where
$U$ denotes the universe of discourse, $x\in U$ and $\mu_{U}(x)$
is a membership function such that $\mu_{U}(x)\in[0,1]$. 
\end{defn}
Hence, the membership function $\mu_{U}(x)$ assigns to every member
$x$ of the set $U$ its grade of membership in $U$.

\textit{Support} of a fuzzy set $\tilde{U}$ is a subset of $U$ for
which the value of the membership function is greater than 0: $Support(\tilde{U})=\lbrace x\in U\vert\mu_{U}(x)>0\rbrace$.

The $\alpha$\textit{-cut} of a fuzzy set $\tilde{U}$ consists of
all elements $x$ in $U$ for which a value of the membership function
is equal to or greater than $\alpha$: $\tilde{U}_{\alpha}=\lbrace x\in U\vert\mu_{U}(x)\geq\alpha\rbrace$ 
\begin{defn}
A \textit{fuzzy number} $\tilde{U}\equiv(U,\mu_{U}(x))$ is a fuzzy
set that is piecewise continuous, achieves the membership value 1
at least once and all its alpha cuts are convex sets \citep{Dubois1980fsas}. 
\end{defn}
\textit{Triangular fuzzy number} (TFN) constitutes a special case
of a fuzzy number such that the membership function has the shape
of a triangle and achieves the value of 1 at exactly one point, see
Figure 1.

Formally, a TFN can be defined as a triplet of real numbers $(l,m,u)$,
$l\leq m\leq u$, where $l$ is the lower bound of the TFN, $m$ is
the value of $x$ for which $\mu_{U}(x)=1$ and $u$ is the upper
bound of the TFN. The membership function of a TFN is given as follows:

\begin{equation}
\mu_{U}(x)=\begin{cases}
0, & \text{if}\ x<l\\
\frac{x-l}{m-l}, & \text{if}\ l\leq x\leq m\\
\frac{u-x}{u-m}, & \text{if}\ m\leq x\leq u\\
0, & \text{if}\ x>u
\end{cases}
\end{equation}

Obviously, the support of a triangular fuzzy number is a closed interval
$[l,u]$.

Arithmetic operations of addition and subtraction with TFNs given
by the Extension principle are closed, which means the result of these
operations is a TFN again. However, multiplication and division of
TFNs do not preserve the triangular shape of a membership function
in general \citep{Brunelli2017aiia}. In practice, piece-wise linear
approximations are used for the two latter operations, see Table 1.

Fuzzy weights need to be normalized in order to generate a unique priority vector (vector of weights). Otherwise, there is an infinite number of fuzzy weights that can be derived from a fuzzy pairwise comparisons matrix, thus making the weights incomparable and impossible to aggregate in a hierarchical structure (if there is any).

\begin{defn}	
Let $ \{\tilde{w}_{i}\,:\, i = 1,...,n \}$ be the set of fuzzy weights. Let $ [l_{i}^{*}(\alpha),u_{i}^{*}(\alpha)] $, $ \alpha \in [0,1] $ represent the $ \alpha $-cut of $ \tilde{w}_{i} $. Then the set $ \{\tilde{w}_{i}\}$ is said to be fuzzy normalized if the following condition holds for every $ \alpha $:
\begin{equation}
[\sum_{i=1}^{n}l_{i}^{*}(\alpha)]\cdot [\sum_{i=1}^{n}u_{i}^{*}(\alpha)] = 1
\end{equation}

\end{defn}
Further on, if the condition above is satisfied only for $\alpha=0$
and $\alpha=1$, then it is called \textit{relaxed normalization}
\citep{Chang1995teon}. In this particular case, we obtain:

\begin{equation}
\left\{ \begin{aligned}\sum_{i=1}^{n}l_{i}\cdot\sum_{i=1}^{n}u_{i}=1(\alpha=0)\\
\sum_{i=1}^{n}m_{i}=1(\alpha=1)
\end{aligned}
\right.
\label{eq:no9}
\end{equation}

\begin{defn}
A \textit{fuzziness} $k(\tilde{U})$ of a fuzzy number $\tilde{U}$
is defined as the area under the membership function $\mu_{U}(x)$
\citep{Kaufmann1988fmmi}:

\begin{equation}
k(\tilde{U})=\int_{support(\tilde{U})}\mu_{U}(x)dx
\end{equation}
\end{defn}
In the case of a triangular fuzzy number $\tilde{U}=(l,m,u)$, the
formula above simplifies to

\begin{equation}
k(\tilde{U})=\frac{u-l}{2}
\end{equation}

In a sense, fuzziness determines the 'specificity' contained in a
given TFN. For example, $\tilde{a}_{ij}=(2,3,4)$ expresses that an
object $i$ is probably two to four times more preferred than an object
$j$, while $\tilde{b}_{ij}=(2,10,90)$ is less specific, since it
says an object $i$ is probably two to ninety times more preferred
than an object $j$. Hence, too much fuzziness (uncertainty) might
be counterproductive, and this notion appears in the formulation of
the mathematical model in Section 4.

If a decision maker wishes to obtain a crisp (real) number from a
triangular fuzzy number, defuzzification must be performed. In this paper we use the following formula called the center of gravity: 

\begin{equation}
c(\tilde{U})=\frac{l+m+u}{3},
\label{eq:center-of-gravity}
\end{equation}
where $c(\tilde{U})$ denotes a crisp value of a triangular fuzzy
number $\tilde{U}=(l,m,u)$.

\begin{figure}[h]
\label{Fig:1}
\begin{centering}
\includegraphics[width=0.8\textwidth]{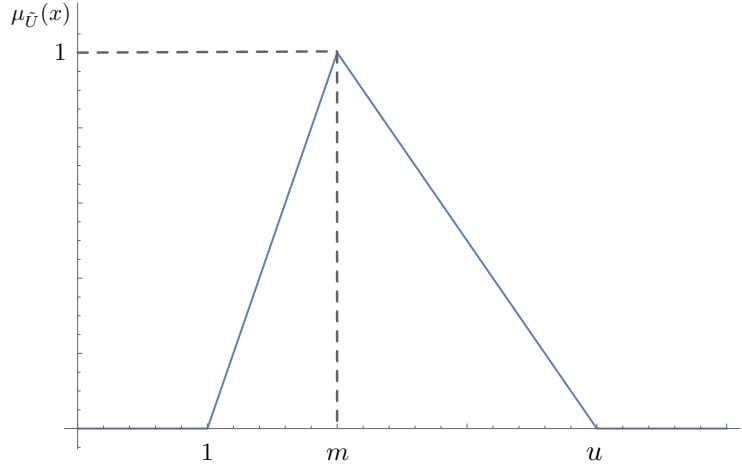}
\par\end{centering}
\caption{A typical membership function of a triangular fuzzy number.}
\end{figure}

\begin{table}[h]
\caption{Standard approximations of arithmetic operations on positive triangular
fuzzy numbers.}
\label{tab:standard-approx} \centering{}%
\begin{tabular}{|c|c|}
\hline 
Arithmetic operation & Definition\tabularnewline
\hline 
$\tilde{A}\oplus\tilde{B}$ & $\langle l_{A}+l_{B},m_{A}+m_{B},u_{A}+u_{B}\rangle$\tabularnewline
\hline 
$\tilde{A}\ominus\tilde{B}$ & $\langle l_{A}-u_{B},m_{A}-m_{B},u_{A}-l_{B}\rangle$\tabularnewline
\hline 
$\tilde{A}\otimes\tilde{B}$ & $\langle l_{A}\cdot l_{B},m_{A}\cdot m_{B},u_{A}\cdot u_{B}\rangle$\tabularnewline
\hline 
$\frac{\tilde{A}}{\tilde{B}}$ & $\langle\frac{l_{A}}{u_{B}},\frac{m_{A}}{m_{B}},\frac{u_{A}}{l_{B}}\rangle$\tabularnewline
\hline 
$\tilde{A}^{-1}$ & $\langle\dfrac{1}{u_{A}},\dfrac{1}{m_{A}},\dfrac{1}{l_{A}}\rangle$\tabularnewline
\hline 
\end{tabular}
\end{table}


\section{Derivation of a priority vector via spanning trees with crisp preferences and confidence levels}

The input for the proposed method consists of an $ n \times n $ (incomplete) \textit{preference}
PC matrix $A=(a_{ij})$, where $a_{ij}\in\mathbb{R}^{+}$, and an $ n \times n $  \textit{confidence}
matrix $C=(c_{ij})$, $c_{ij}\in R_{0}^{+}$, where the values of
$c_{ij}$ express the confidence of corresponding values $a_{ij}$.
It is assumed that both $a_{ij}$ and $c_{ij}$ take values from some
suitable comparison scale set in advance (for instance, one can apply
Saaty's scale from 1 to 9 for $a_{ij}$ values and a 5-item ordinal
linguistic scale for the values of $c_{ij}$).

The method starts with the enumeration of all spanning trees from
a graph representation of the matrix $A$ by a suitable algorithm,
see e.g. \citep{Chakraborty2019afga}. For each spanning tree, a unique (normalized) priority vector $w$ is calculated from $(n-1)$ graph edges ($a_{ij}$ values) and the
normality condition $\sum_{i=1}^{n}w_{i}=1$. Subsequently, this priority vector is assigned its absolute and relative reliability according to relations (5\textendash 6). 
At last, the final priority vector $w^{(f)}$ is obtained as a weighted
arithmetic mean of the priority vectors $w^{(STi)}$ obtained from $ N $
individual spanning trees, where the relative reliability $R^{(i)}$
of the i-th spanning tree plays the role of weights:
\begin{equation}
w^{(f)}=\sum_{i=1}^{N}R^{(i)}\cdot w^{(STi)}
\end{equation}
Based on the values of $w_{i}^{(f)}$ all compared objects can be
ranked (ordered from the best to the worst).

To summarize the whole process, each spanning tree of the graph representation of the preference matrix leads to a corresponding (individual) priority vector, and the final priority vector of all compared objects is obtained as an aggregation of priority vectors over all spanning trees. In addition, a confidence value is assigned to each pairwise comparison. Consequently, each priority vector from a given spanning tree has a reliability value, which is the product of the confidence levels of individual pairwise comparisons forming the tree. In this way, each priority vector acquires its weight equal to the reliability in the final aggregation process, where the greater reliability means the greater weight. The method with the fuzzy extension is analogous; only arithmetic operations over triangular fuzzy numbers are performed differently than for crisp numbers.
\indent The method is illustrated on the following numerical example. 
\begin{example}
Consider the incomplete multiplicative PC matrix $A$ and the confidence
matrix $C$ shown below. While $a_{ij}$, elements of $A$, express
a degree of preference of an object $i$ over an object $j$, the
elements $c_{ij}$ of a matrix $C$ express the confidence of an expert
of the correct value of $a_{ij}$. The sign "{*}" denotes a missing
element. Notice that while $a_{ij}=1/a_{ji},\forall i,j$ (the reciprocity
condition), we have $c_{ij}=c_{ji},\forall i,j$ (the symmetry condition).

\[
\mathbf{A}=\left[\begin{array}{cccc}
1 & 4 & 1/2 & *\\
1/4 & 1 & 2 & *\\
2 & 1/2 & 1 & 5\\
* & * & 1/5 & 1
\end{array}\right],
\]

\[
\mathbf{C}=\left[\begin{array}{cccc}
- & 4 & 3 & *\\
4 & - & 1 & *\\
3 & 1 & - & 2\\
* & * & 2 & -
\end{array}\right],
\]

The graph representation of the matrix $A$ is shown in Figure 2,
and all three spanning trees are shown in Figure \ref{Fig:3}. Each
spanning tree provides a unique vector of weights $w=(w_{1},...,w_{4})$,
$\Vert w\Vert=1$, derived from ratios of all $a_{ij}$ included in
a given spanning tree.

In the case of the spanning tree number 1 (see Figure \ref{Fig:3}),
we obtain the following set of equations:

\begin{equation}
\left\{ \begin{aligned}w_{1}=4w_{2},\\
w_{2}=2w_{3},\\
w_{4}=5w_{3},\\
\sum_{i=1}^{4}w_{i}=1
\end{aligned}
\right.
\end{equation}

The solution of the system of equations (14) is: $w^{(ST1)}=(0.5,0.125,0.0625,0.3125)$,
see also Table \ref{tab:priority-vectors}, second row.

The absolute reliability of this priority vector is estimated as follows:

$r^{(1)}=c_{12}\cdot c_{23}\cdot c_{34}=4\cdot1\cdot2=8$,

And the relative reliability is given as:

$R^{(1)}=8/(8+24+6)=8/38=0.211$.

Further on, $R^{(2)}=24/38=0.632$ and $R^{(3)}=6/38=0.158$

All priority vectors and their relative reliability are shown in Table
\ref{tab:priority-vectors}. As can be seen, more reliable priority
vectors (vectors from spanning trees with higher reliability) achieve
higher weights in the final aggregation.

\begin{table}[h]
\caption{Priority vectors and their relative reliability for all three spanning
trees}
\label{tab:priority-vectors} \centering{}%
\begin{tabular}{|c|c|c|}
\hline 
Spanning tree & Priority vector & $(R)$\tabularnewline
\hline 
1 & $w^{(ST1})$ = (0.5, 0.125, 0.0625, 0.3125) & 0.211\tabularnewline
\hline 
2 & $w^{(ST2})$ = (0.274, 0.068, 0.548, 0.110) & 0.632\tabularnewline
\hline 
3 & $w^{(ST3})$ = (0.135, 0.541, 0.270, 0.054) & 0.158\tabularnewline
\hline 
\end{tabular}
\end{table}

The final priority vector is estimated as a weighted arithmetic mean
(13), where the reliability plays the role of weights of each spanning
tree priority vector:

$w^{(f)}=\sum_{i=1}^{3}R^{(i)}\cdot w^{(STi)}=0.211\cdot(0.5,0.125,0.0625,0.3125)$
+ $0.632\cdot(0.274,0.068,0.548,0.110)+0.158\cdot(0.135,0.541,0.270,0.054)=(0.300,0.155,0.402,0.144)$.

According to the final priority vector, the ranking of compared objects
is as follows: (3,1,2,4).

\selectlanguage{english}%
\begin{figure}[h]
\label{Fig:1-1}
\begin{centering}
\includegraphics{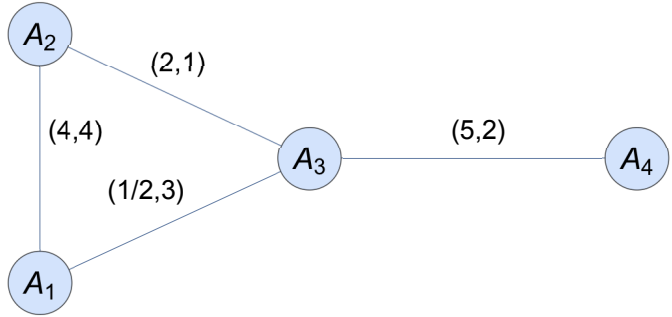}
\par\end{centering}
\caption{Graph representation\protect\footnotemark of a PC matrix $A$.}
\end{figure}

\selectlanguage{american}%
\footnotetext{We assume that the edge labels correspond to the increasing
order of vertices. In other words, the label between $A_{i}$ and
$A_{j}$ is set to $a_{ij}$ if $i\leq j$, to $a_{ji}$ otherwise.
}

\begin{figure}[h]
\label{Fig:2}
\begin{centering}
\subfloat[ST1]{\begin{centering}
\includegraphics{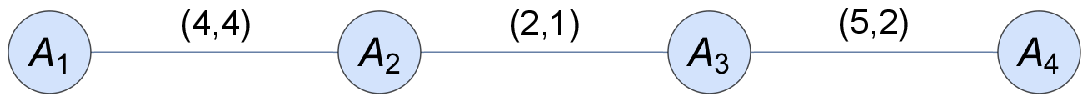}
\par\end{centering}
}
\par\end{centering}
\begin{centering}
\subfloat[ST2]{\begin{centering}
\includegraphics{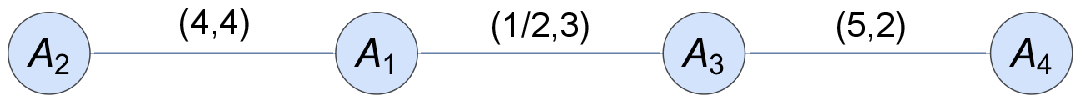}
\par\end{centering}
}
\par\end{centering}
\begin{centering}
\subfloat[ST3]{\begin{centering}
\includegraphics{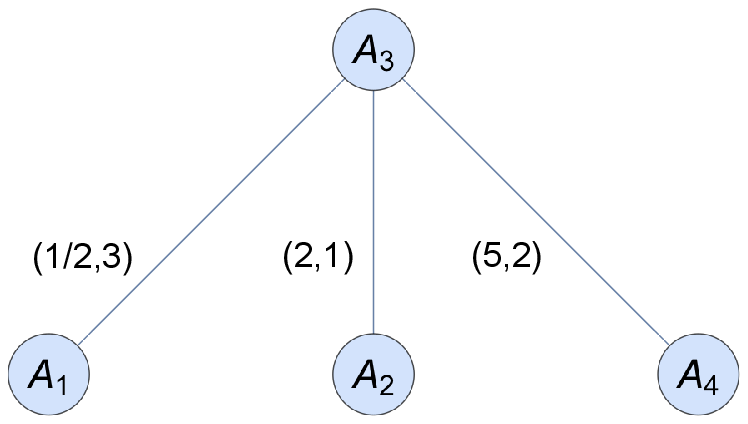}
\par\end{centering}
}
\par\end{centering}
\caption{Three spanning trees of $A$}
\label{Fig:3}
\end{figure}
\end{example}

A comparison of the proposed method with Harker's method \citep{Harker1987amoq}
and the geometric mean method modified for incomplete PC matrices
\citep{kulakowski2020otgm} is provided below, but without confidence
levels, of course. 
\begin{example}
. Consider the same PC matrix as in Example 15. Harker's method starts
with the construction of an auxiliary matrix $B=[b_{ij}]$ associated
with $A$, such that

\[
b_{ij}=\begin{cases}
0 & \text{if}\,\,a_{ij}=?\,\,\text{and}\,\,i\neq j\\
c_{ij} & \text{if}\,\,a_{ij}\neq?\,\,\text{and}\,\,i\neq j\\
s_{i}+1 & \text{if}\,\,i=j
\end{cases},
\]

where $s_{i}$ denotes the number of missing entries in the i-th row
of $A$, given as: 
\[
B=\left(\begin{array}{cccc}
2 & 4 & \frac{1}{2} & 0\\
\frac{1}{4} & 2 & 2 & 0\\
2 & \frac{1}{2} & 1 & 5\\
0 & 0 & \frac{1}{5} & 3
\end{array}\right).
\]

By solving the equation 
\[
Bw=\lambda_{\textit{max}}w
\]

we obtain the ranking vector $w_{\textit{hm}}$, which after the appropriate
rescaling equals

\[
w_{\textit{hm}}=\left(\begin{array}{c}
0.416\\
0.252\\
0.298\\
0.0335
\end{array}\right).
\]

Similarly, using the geometric mean method for incomplete PC matrices
\citep{kulakowski2020otgm} we obtain

\[
w_{\textit{gm}}=\left(\begin{array}{c}
\frac{e^{\frac{1}{12}(4\ln(2)+3\ln(5))}}{\sqrt[4]{5}+\frac{1}{5^{3/4}}+e^{\frac{1}{12}(3\ln(5)-4\ln(2))}+e^{\frac{1}{12}(4\ln(2)+3\ln(5))}}\\
\frac{e^{\frac{1}{12}(3\ln(5)-4\ln(2))}}{\sqrt[4]{5}+\frac{1}{5^{3/4}}+e^{\frac{1}{12}(3\ln(5)-4\ln(2))}+e^{\frac{1}{12}(4\ln(2)+3\ln(5))}}\\
\frac{\sqrt[4]{5}}{\sqrt[4]{5}+\frac{1}{5^{3/4}}+e^{\frac{1}{12}(3\ln(5)-4\ln(2))}+e^{\frac{1}{12}(4\ln(2)+3\ln(5))}}\\
\frac{1}{5^{3/4}\left(\sqrt[4]{5}+\frac{1}{5^{3/4}}+e^{\frac{1}{12}(3\ln(5)-4\ln(2))}+e^{\frac{1}{12}(4\ln(2)+3\ln(5))}\right)}
\end{array}\right),
\]

i.e.

\[
w_{\textit{gm}}=\left(\begin{array}{c}
\frac{5\cdot2^{2/3}}{5+6\sqrt[3]{2}+5\cdot2^{2/3}}\\
\frac{5}{5+6\sqrt[3]{2}+5\cdot2^{2/3}}\\
\frac{5\sqrt[3]{2}}{5+6\sqrt[3]{2}+5\cdot2^{2/3}}\\
\frac{\sqrt[3]{2}}{5+6\sqrt[3]{2}+5\cdot2^{2/3}}
\end{array}\right),
\]

and finally 
\[
w_{\textit{gm}}=\left(\begin{array}{c}
0.387\\
0.2439\\
0.307\\
0.061
\end{array}\right).
\]
\end{example}
As can be seen, both Harker's method and the modified geometric mean
method yield the same ranking of the compared objects: (1,3,2,4),
which slightly differs from the ranking (3,1,2,4) obtained from the
model with confidence levels. This difference can be attributed to
the application of additional information contained in the confidence
matrix of the proposed model.

\section{Derivation of a priority vector via spanning trees with fuzzy preferences and confidence levels}

The input for the proposed method consists of a (incomplete) $n\times n$
preference PC matrix $\tilde{A}=(\tilde{a}_{ij})$, where $\tilde{a}_{ij}$
has the form of triangular fuzzy numbers, and an $ n \times n $ confidence matrix
$\tilde{C}=(\tilde{c}_{ij})$, where $\tilde{c}_{ij}$ are triangular
fuzzy numbers as well. In addition, it is assumed that (fuzzy) scales
for both pairwise comparisons and confidence levels were selected
in advance.

The method begins with the enumeration of all spanning trees from
a graphical representation of preferences provided in a PC matrix
in the same way as in the previous section.

Let $\tilde{w}_{i}\equiv(l_{i},m_{i},u_{i})$, $i\in\lbrace1,...,n\rbrace$
denote weights (coordinates of a priority vector) in the form of TFNs
of all compared objects with respect to a given spanning tree that
is formed by $(n-1)$ edges denoted as $[k,l]$, $k,l\in\lbrace1,...,n\rbrace$.


Then, for each spanning tree the following mathematical model for weights derivation comprised of four parts (I-IV) is constructed:

\begin{equation}
(I):\tilde{w}_{i}=\tilde{a}_{ij}\cdot\tilde{w}_{j},j=i+1,i\in\lbrace1,...,n\rbrace
\end{equation}

\begin{equation}
(II):0<l_{i}\leq m_{i}\leq u_{i},i\in\lbrace1,...,n\rbrace
\end{equation}

\begin{equation}
\left\{ \begin{aligned}(III):\\
(\sum_{i=1}^{n}l_{i})\cdot(\sum_{i=1}^{n}u_{i})=1\\
\sum_{i=1}^{n}m_{i}=1
\end{aligned}
\right.
\label{eq:no17}
\end{equation}

\begin{equation}
\left\{ \begin{aligned}(IV):\\
min\sum_{i=1}^{n}(u_{i}-l_{i})
\end{aligned}
\right.
\end{equation}

Part (I) is the `Preference part' of the model. It consists of $(n-1)$
fuzzy equations that can be decomposed into $3(n-1)$ crisp equations
for the $3n$ unknown values of $(l_{i},m_{i},u_{i})$. Since there
are $3n$ unknown variables and only $3(n-1)$ equations, there are
3 degrees of freedom.

Part (II) is the `Structural part' of the model. It places constraints
on the values of lower bounds, middle values and upper bound of all
triangular fuzzy numbers. If this part was missing, a solution that
contradicts $0<l_{i}\leq m_{i}\leq u_{i}$ inequalities could be obtained
(this was noted as early as the Chang and Lee paper from 1995,\citep{Chang1995teon}
).


Part (III) is the `Normalization part' of the model since it is usually desired that weights 
(a priority vector) are normalized. Equations (\ref{eq:no17}) are identical to 
Equations (\ref{eq:no9}) expressing the condition of the relaxed normalization. 
These two additional equations reduce the degrees of freedom of the model to 1.	

Finally, Part (IV) has the following purpose: the model has one degree
of freedom, hence one can expect an infinite number of solutions (if
the set of inequalities in Part (II), which provides other constraints
for the solution, is satisfied) depending on one parameter (it can
be $u_{1}$, for example). If this happens, then condition (IV) provides
a unique solution with minimal 'fuzziness' (a solution in the form
of a 'narrow' triangular fuzzy number provides more specific information
than a 'wide' triangular fuzzy number).

After a priority vector $\tilde{w}$ from a spanning tree $ST$ with the set of edges $E^{(ST)}$ is acquired, its absolute reliability $\tilde{r}$ is given as follows:

\begin{equation}
\tilde{r}(ST)=\lbrace\prod\tilde{c}_{ij}\vert[ij]\in E^{(ST)}\rbrace
\end{equation}

The corresponding relative reliability $\tilde{R}$ of a priority vector $\tilde{w}$ from a spanning tree $ST$  is defined as:

\begin{equation}
\tilde{R}(ST)=\frac{\tilde{r}(ST)}{\sum_{i=1}^{N}\tilde{r}(STi)},
\end{equation}

where $N$ is the number of spanning trees.

Note that operations of multiplication and division are fuzzy operations
defined in Table \ref{tab:standard-approx}.

At last, the final (fuzzy) priority vector $\tilde{w}^{(f)}$ is obtained
as a weighted arithmetic mean of priority vectors $\tilde{w}^{(STi)}$
obtained from $ N $ individual spanning trees $(STi)$, where the relative
reliability $\tilde{R}^{(i)}$ of the i-th spanning tree plays the
role of weights:

\begin{equation}
\tilde{w}^{(f)}=\sum_{i=1}^{N}\tilde{R}^{(i)}\otimes\tilde{w}^{(ST_{i})},
\end{equation}
where '$\otimes$' denotes fuzzy multiplication.

The use of the method is illustrated in the following example. 
\begin{example}
Consider the incomplete multiplicative matrix PC matrix $\tilde{A}=(\tilde{a}_{ij})$,
as shown below, and its corresponding confidence matrix $\tilde{C}=(\tilde{c}_{ij})$,
where the elements of both $\tilde{A}$ and $\tilde{C}$ are in the
form of triangular fuzzy numbers. In order to enable a comparison
of methods applied in Examples 15 and 16, we only transformed the crisp
values from Example 15 to the corresponding triangular fuzzy numbers.
We applied Saaty's 9-point scale from \citep{Ayhan2013afaa}, see
Fig. 4, and the 5-item fuzzy DEMATEL linguistic scale from \citep{Kurniawati2015pios},
see Table \ref{tab:the-fuzzy-linguistic}.

\bigskip{}

From the the first spanning tree (see Figure \ref{Fig:3}), we obtain
the following mathematical model for the priority vector $\tilde{w}=(\tilde{w}_{1},...,\tilde{w}_{n})$:

\begin{equation}
\left\{ \begin{aligned}\tilde{w}_{1}=\tilde{a}_{12}\otimes\tilde{w}_{2}=(3,4,5)\otimes\tilde{w}_{2}\\
\tilde{w}_{2}=\tilde{a}_{23}\otimes\tilde{w}_{3}=(1,2,3)\otimes\tilde{w}_{3}\\
\tilde{w}_{3}=\tilde{a}_{34}\otimes\tilde{w}_{4}=(4,5,6)\otimes\tilde{w}_{4}
\end{aligned}
\right.
\end{equation}

\begin{equation}
\left\{ \begin{aligned}0<l_{1}\leq m_{1}\leq u_{1}\\
0<l_{2}\leq m_{2}\leq u_{2}\\
0<l_{3}\leq m_{3}\leq u_{3}\\
0<l_{4}\leq m_{4}\leq u_{4}
\end{aligned}
\right.
\end{equation}

\begin{equation}
\left\{ \begin{aligned}(l_{1}+l_{2}+l_{3}+l_{4})\cdot(u_{1}+u_{2}+u_{3}+u_{4})=1\\
m_{1}+m_{2}+m_{3}+m_{4}=1
\end{aligned}
\right.
\end{equation}

\begin{equation}
min[(u_{1}+u_{2}+u_{3}+u_{4})-(l_{1}+l_{2}+l_{3}+l_{4})]
\end{equation}

The solution of the model (obtained by \textit{Mathematica Solver})
is as follows:

\begin{equation}
\left\{ \begin{aligned}\tilde{w}_{1}=(0.214,0.714,2.09)\\
\tilde{w}_{2}=(0.071,0.179,0.417)\\
\tilde{w}_{3}=(0.071,0.089,0.089)\\
\tilde{w}_{4}=(0.018,0.018,0.023)
\end{aligned}
\right.
\end{equation}

The absolute and relative reliabilities of this priority vector are:

$\tilde{r}(1)=\tilde{c}_{12}\otimes\tilde{c}_{23}\otimes\tilde{c}_{34}=(0.75,1,1)\otimes(0,0.25,0.5)\otimes(0.25,0.5,0.75)=(0,0.125,0.375)$, and $\tilde{R}(1)=(0,0.211,4)$.

The results for all three spanning trees are summarized in Table \ref{Outcomes}
along with the final priority vector $\tilde{w}^{(f)}$ obtained from
relation (21). Here, only the calculation of the first coordinate
of $\tilde{w}^{(f)}$ is shown (the rest is analogous):

\[ 
\begin{split}
\tilde{w}_{1}^{(f)}=&(0,0,211,4) \otimes(0.214,0.714,2.09) \oplus(0.063,0.632,8)\otimes(0.206,0.274,0.401)\\ & \oplus(0,0.158,4)  \otimes(0.135,0.135,0.169)\\&=(0.013,0.345,12.23) 
\end{split}
 \]

Once the final fuzzy weights (final priority vector) $ {\tilde{w}_{i}^{(f)} }$ are obtained, see the lower part of Table \ref{Outcomes}, the defuzzification by relation (\ref{eq:center-of-gravity}) and subsequent normalization is performed, see the bottom of Table \ref{Outcomes} for results.
Now, crisp and normalized weights of each object can be easily compared and objects can be ranked for the best to the worst. It should be noted that the defuzzification is not necessary since a comparison of fuzzy numbers is also possible, but it is far from being straightforward and thus less suitable for practical applications.

After the defuzzification by relation (\ref{eq:center-of-gravity}), the final (crisp) priority vector is: $w=(0.353$, $0.240$, $0.348$, $0.059)$. The ranking of objects (from the best
to the worst) is as follows:(1,3,2,4), which is the same as the ranking
obtained from Harker's method and the GMM method in Example 16, but
different from the result obtained via the spanning trees method with
crisp confidence levels in Example 15.

Notice that from the table of arithmetic operations, it follows that
$\tilde{A}\otimes\frac{1}{\tilde{A}}\neq\tilde{1}$. Therefore, in
the system of equations (15) and (22), respectively, one should consistently
multiply TFNs so that all $\tilde{c}_{ij}$ have the lower bound $l_{ij}\geq1$
(or vice versa) to obtain replicable results.

\[
\mathbf{A}=\left[\begin{array}{cccc}
(1,1,1) & (3,4,5) & (1/3,1/2,1) & *\\
(1/5,1/4,1/3) & (1,1,1) & (1,2,3) & *\\
(1,2,3) & (1/3,1/2,1) & (1,1,1) & (4,5,6)\\
* & * & (1/6,1/5,1/4) & (1,1,1)
\end{array}\right],
\]

\[
\mathbf{C}=\left[\begin{array}{cccc}
- & (0.75,1,1) & (0.5,0.75,1) & *\\
(0.75,1,1) & - & (0,0.25,0.5) & *\\
(0.5,0.75,1) & (0,0.25,0.5) & - & (0.25,0.5,0.75)\\
* & * & (0.25,0.5,0.75) & (1,1,1)
\end{array}\right],
\]

\begin{figure}[h]
\begin{centering}
\includegraphics[scale=0.3,angle=270]{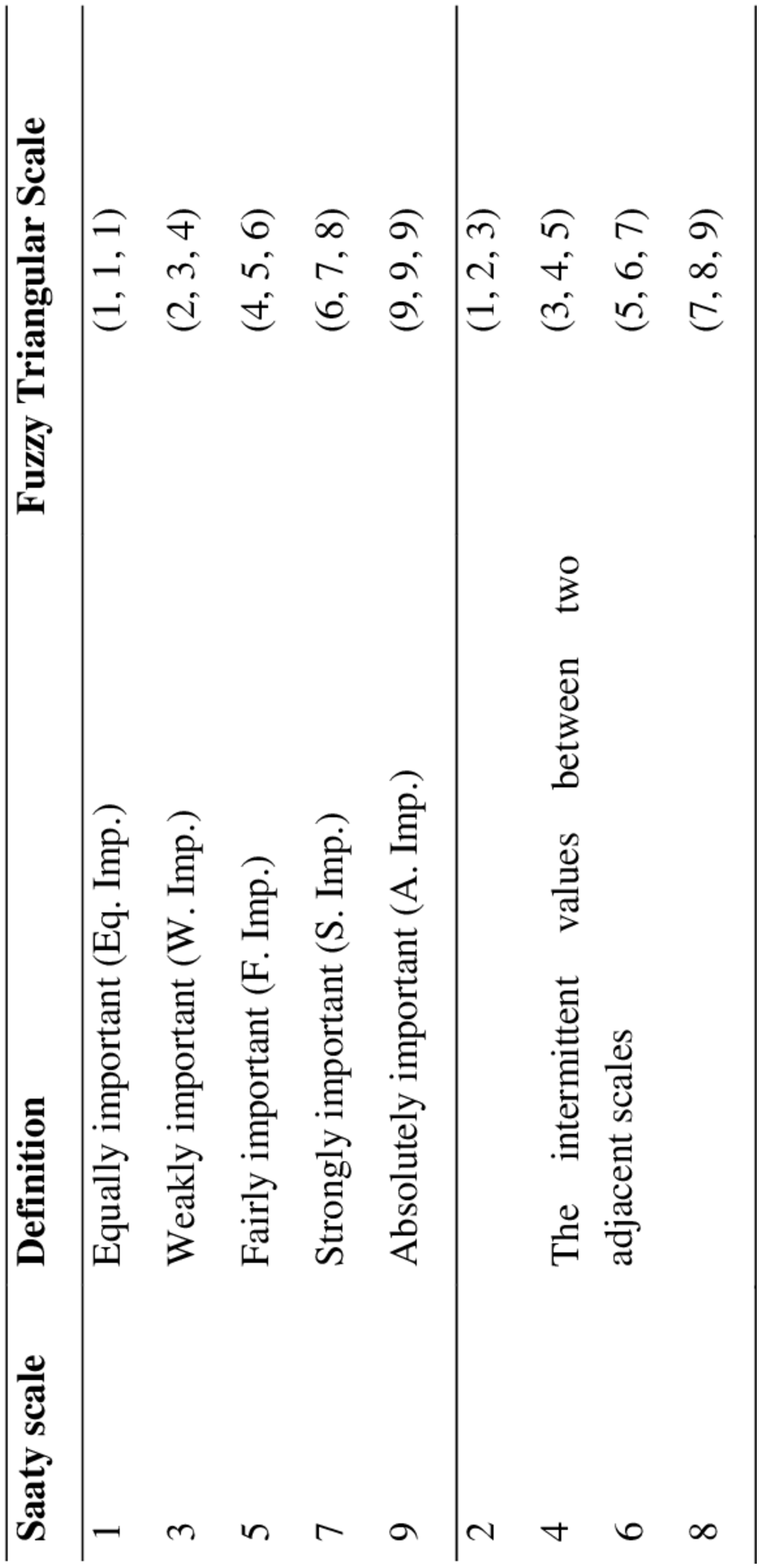}
\par\end{centering}
\caption{Saaty's fuzzy 1-9 scale. Source: \citep{Ayhan2013afaa}.}
\label{Fig:1-2}
\end{figure}

\begin{table}[h]
\caption{The fuzzy linguistic scale. Source: \citep{Kurniawati2015pios}}
\label{tab:the-fuzzy-linguistic} \centering{}%
\begin{tabular}{|c|c|c|}
\hline 
Linguistic term & Crisp value & TFN\tabularnewline
\hline 
zero confidence & 0 & (0,0,0.25)\tabularnewline
\hline 
low confidence & 1 & (0,0.25,0.5)\tabularnewline
\hline 
moderate confidence & 2 & (0.25, 0.5, 0.75)\tabularnewline
\hline 
strong confidence & 3 & (0.5, 0.75, 1)\tabularnewline
\hline 
absolute confidence & 4 & (0.75, 1, 1)\tabularnewline
\hline 
\end{tabular}
\end{table}
\end{example}
\begin{table}[!h]
\caption{Example 17 - The outcomes}
\label{Outcomes} \centering{}%
\begin{tabular}{|c|c|c|c|}
\hline 
Spanning tree 1 & l & m & u\tabularnewline
\hline 
$\tilde{w}_{1}$ & 0.2143 & 0.7143 & 2.087\tabularnewline
$\tilde{w}_{2}$ & 0.0714 & 0.1786 & 0.4174\tabularnewline
$\tilde{w}_{3}$ & 0.0714 & 0.0893 & 0.0893\tabularnewline
$\tilde{w}_{4}$ & 0.0179 & 0.0179 & 0.0232\tabularnewline
$\tilde{R}$ & 0 & 0.2105 & 4\tabularnewline
\hline 
Spanning tree 2 & l & m & u\tabularnewline
\hline 
$\tilde{w}_{1}$ & 0.2055 & 0.2740 & 0.4008\tabularnewline
$\tilde{w}_{2}$ & 0.0685 & 0.0685 & 0.0802\tabularnewline
$\tilde{w}_{3}$ & 0.2055 & 0.5479 & 0.2025\tabularnewline
$\tilde{w}_{4}$ & 0.0514 & 0.1096 & 0.2004\tabularnewline
$\tilde{R}$ & 0.0625 & 0.6315 & 8\tabularnewline
\hline 
Spanning tree 3 & l & m & u\tabularnewline
\hline 
$\tilde{w}_{1}$ & 0.1351 & 0.1351 & 0.1687\tabularnewline
$\tilde{w}_{2}$ & 0.1351 & 0.5405 & 1.518\tabularnewline
$\tilde{w}_{3}$ & 0.1351 & 0.2703 & 0.5060\tabularnewline
$\tilde{w}_{4}$ & 0.0338 & 0.0541 & 0.0843\tabularnewline
$\tilde{R}$ & 0 & 0.1578 & 4\tabularnewline
\hline 
Final priority vector & l & m & u\tabularnewline
\hline 
$\tilde{w}^{(f)}_{1}$ & 0.0128 & 0.3448 & 12.2292\tabularnewline
$\tilde{w}^{(f)}_{2}$ & 0.0043 & 0.1662 & 8.3832\tabularnewline
$\tilde{w}^{(f)}_{3}$ & 0.0128 & 0.4075 & 12.0012\tabularnewline
$\tilde{w}^{(f)}_{4}$ & 0.0032 & 0.0815 & 2.0332\tabularnewline
\hline 
Final crisp priority vector & after defuzzification & after normalization & rank\tabularnewline
\hline 
$w_{1}$ & 4.195 & 0.353 & 1\tabularnewline
$w_{2}$ & 2.851 & 0.240 & 3\tabularnewline
$w_{3}$ & 4.141 & 0.348 & 2\tabularnewline
$w_{4}$ & 0.706 & 0.059 & 4\tabularnewline
\hline 

\end{tabular}
\end{table}



\section{Discussion}

\label{sec12} 
In the proposed approach, we use the term 'confidence' to estimate 
the belief in the accuracy of individual pairwise comparisons provided 
by an expert, which allows an expert to assess the quality of their judgments. 
This feature remains neglected in the standard pairwise comparisons methods 
such as the analytic hierarchy/network process, where it is assumed all pairwise 
comparisons have the same confidence. However, in real-world problems involving 
(not only) brand new complex phenomena or emerging technologies, it is natural 
that even experts with the best knowledge might be, in some situations, 
less confident about their judgments. Moreover, our approach allows straightforward 
extensions to a multiple-criteria single or group decision-making, where experts' 
assessments concerning different criteria may be performed with individually 
selected confidence.

Apart from our method's already mentioned advantages
(there is no need to fill the incomplete PC matrix and confidence
levels enable modeling uncertainty), the method has some limitations,
specifically it requires that all spanning trees are found (there
are a lot of algorithms \citep{Gross2018gtaI} that can do the task)
and a priority vector and its reliability estimated for each spanning
tree. However, the number of spanning trees, especially for complete
or almost complete PC matrices (graphs) could be very high, rendering
the method excessively time-consuming and computationally demanding
to use. That's why we propose using the method primarily for incomplete
matrices.

Before a decision maker applies our method, we suggest approximately
estimating the number of spanning trees that would be necessary to
process. For this, we postulate the following simple combinatorial
theorem which gives an upper bound for the number of spanning trees
$t(A)$ associated with a PC matrix $A$ which is based only on the
number of vertices, the number of edges and the number of vertices
with the degree equal to 1 and does not require formation of the Laplacian
matrix and calculation of its cofactors or eigenvalues. 
\begin{thm}
Let $A$ be a PC matrix, let $G$ be a graph representation of $A$
and let $t(G)$ be the number of spanning trees of $G$. Further on,
let \emph{$n,e,k$} denote the number of vertices, the number of edges
and the number of vertices with the degree equal to 1 (called leaf
vertex or end vertex, since it is connected to the rest of a graph
with only one edge), respectively, of the graph $G$. Then: 
\begin{equation}
t(G)\leq\binom{e-k}{n-k-1}
\end{equation}
\end{thm}
\begin{proof}
A graph with $n$ vertices and $e$ edges requires each spanning tree
to have $n-1$ edges. There are $\binom{e}{n-1}$ possible ways to
perform the selection. Of course, the number of spanning trees would
be the same or lower than $\binom{e}{n-1}$ since some choices of
$n-1$ edges may form circles, leaving the rest of the graph unconnected.
Therefore, $t(G)\leq\binom{e}{n-1}$. Further on, if there is a vertex
with the degree equal to 1, the edges from this vertex have to be
included in each spanning tree. Thus, the combinatorial number $\binom{e}{n-1}$
describing a selection of $n-1$ edges from $e$ changes to $\binom{e-1}{n-2}$,
since both the number of edges available and the number of edges necessary
to form a spanning tree both decrease by 1. Apparently, when there
are $k$ vertices with the degree of one, the combinatorial number
becomes $\binom{e-k}{n-k-1}$. 
\end{proof}
For a cyclic graph ($e=n,k=0$) we get $t(G)\leq\binom{e-k}{n-k-1}=\binom{n}{n-1}=n$,
hence the equality holds. The same relationship holds if a graph is
a spanning tree itself ($e=n-1,k=2$): $t(G)\leq\binom{e-k}{n-k-1}=\binom{n-1-2}{n-2-1}=1$.

Using Theorem 18, a decision maker can decide beforehand if the proposed
method would be feasible. The following Example 4 illustrates the
use of Theorem 18. 
\begin{example}
. Consider the following incomplete real PC matrix $A$ of the order
$n=7$ with a graph representation depicted in Figure 5. Estimate
the upper bound for the number of spanning trees necessary to perform
the proposed method.

\[
\mathbf{A}=\left[\begin{array}{ccccccc}
1 & 5 & * & * & * & * & 0.25\\
0.2 & 1 & 4 & 0.25 & 4 & * & *\\
* & 0.25 & 1 & 0.5 & * & * & *\\
* & 8 & 2 & 1 & 2 & * & 6\\
* & 0.25 & * & 0.5 & 1 & * & *\\
* & * & * & * & * & 1 & 3\\
4 & * & * & 0.167 & * & 0.333 & 1
\end{array}\right],
\]

From the matrix $A$ it is clear that $n=7$, $e=9$ and $k=1$ (the
object number 6 is connected to the object number 7 only). From Theorem
18 it follows that $t(G)\leq\binom{e-k}{n-k-1}=\binom{8}{5}=56$.

Hence, there are only 56 spanning trees, or less. The Kelmans and
Chelnokov formula \citep{Kelmans1974acpo} for the exact number of
spanning trees gives the value of 28.

\begin{figure}[H]
\begin{centering}
\includegraphics{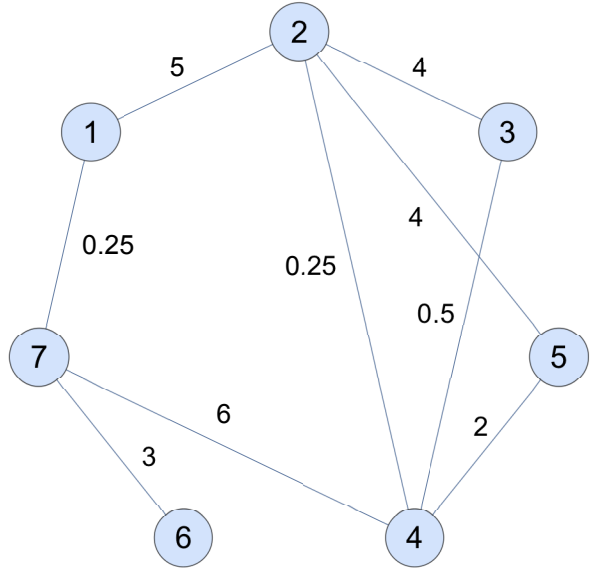}
\par\end{centering}
\caption{Example matrix graph - undirected}
\end{figure}
\end{example}

\section{Conclusion}

\label{sec13}

The aim of the paper was to introduce a novel approach for the derivation
of a priority vector from an incomplete PC matrix when additional
information on the confidence of the pairwise comparisons is provided.
The method has several advantages, namely it does not require the
completion of a PC matrix, thus reflecting original preferences without
the distortion. Secondly, it enables sophisticated modeling of uncertainty
of pairwise comparisons using embedded confidence levels expressed
by crisp (real) numbers or triangular fuzzy numbers.

Our approach opens several interesting research directions that can
be pursued in the future. Firstly, numerical simulations can be carried
out for the comparison of the proposed method and other methods for
the derivation of a priority vector. Secondly, the problem of the
\textit{condition of the order preservation}, see e.g. \citep{BanaeCosta2008acao,Kulakowski2019witc,Mazurek2020sotc}
, or \citep{Mazurek2019snpo} can be investigated in the context of
the proposed method, with the interplay between preferences and confidences
on one hand and inconsistency on the other hand. Further on, future
research may focus on the extension of the proposed method, for instance
towards Z-numbers, and also to the design of a user-friendly and free
online software tool for facilitating the method's computations. 
In addition, we will also want to investigate the relationship between the level of confrontation and the decision security of the model. We hope that the introduced extension can effectively detect manipulation attempts in the pairwise comparison method.

\section*{Acknowledgements}

The paper was supported by the Grant Agency of the Czech Republic no. 21-03085S, 
and by SODA Grant no. 2021/41/B/HS4/03475 founded by the National Science Center Poland (NCN).

\section*{Literature }

\bibliographystyle{plainnat}
\bibliography{papers_biblio_reviewed}

\begin{thebibliography}{55}
\providecommand{\natexlab}[1]{#1}
\providecommand{\url}[1]{\texttt{#1}}
\expandafter\ifx\csname urlstyle\endcsname\relax
  \providecommand{\doi}[1]{doi: #1}\else
  \providecommand{\doi}{doi: \begingroup \urlstyle{rm}\Url}\fi

\bibitem[Alonso et~al.(2008)Alonso, Chiclana, Herrera, Viedma-Herrera,
  Alcal{\'a}-Fdez, and Porcel]{Alonso2008acbp}
S.~Alonso, F.~Chiclana, F.~Herrera, E.~Viedma-Herrera, J.~Alcal{\'a}-Fdez, and
  C.~Porcel.
\newblock {A consistency-based procedure to estimate missing pairwise
  preference values}.
\newblock \emph{International Journal of Intelligent Systems}, 23\penalty0
  (2):\penalty0 155--175, February 2008.

\bibitem[Ayhan(2013)]{Ayhan2013afaa}
M.~B. Ayhan.
\newblock A fuzzy {AHP} approach for supplier selection problem: {A} case study
  in a gear motor company.
\newblock \emph{CoRR}, abs/1311.2886, 2013.
\newblock URL \url{http://arxiv.org/abs/1311.2886}.

\bibitem[Azadeh et~al.(2013)Azadeh, Saberi, Atashbar, Chang, and
  Pazhoheshfar]{Azadeh2013zahp}
A.~Azadeh, M.~Saberi, N.~Z. Atashbar, E.~Chang, and P.~Pazhoheshfar.
\newblock {Z-AHP: A Z-number extension of fuzzy analytical hierarchy process}.
\newblock \emph{IEEE International Conference on Digital Ecosystems and
  Technologies}, pages 141--147, 2013.
\newblock ISSN 21504938.
\newblock \doi{10.1109/DEST.2013.6611344}.

\bibitem[Bana~e Costa and Vansnick(2008)]{BanaeCosta2008acao}
C.~A. Bana~e Costa and J.~Vansnick.
\newblock {A critical analysis of the eigenvalue method used to derive
  priorities in AHP}.
\newblock \emph{European Journal of Operational Research}, 187\penalty0
  (3):\penalty0 1422--1428, June 2008.

\bibitem[Ben{\'{i}}tez et~al.(2019)Ben{\'{i}}tez, Carpitella, Certa, and
  Izquierdo]{Benitez2019moup}
J.~Ben{\'{i}}tez, S.~Carpitella, A.~Certa, and J.~Izquierdo.
\newblock {Management of uncertain pairwise comparisons in AHP through
  probabilistic concepts}.
\newblock \emph{Applied Soft Computing Journal}, 78:\penalty0 274--285, 2019.
\newblock ISSN 15684946.
\newblock \doi{10.1016/j.asoc.2019.02.020}.

\bibitem[Boz{\'o}ki and Tsyganok(2019)]{Bozoki2019tlls}
S.~Boz{\'o}ki and V.~Tsyganok.
\newblock The (logarithmic) least squares optimality of the arithmetic
  (geometric) mean of weight vectors calculated from all spanning trees for
  incomplete additive (multiplicative) pairwise comparison matrices.
\newblock \emph{International Journal of General Systems}, 48\penalty0
  (4):\penalty0 362--381, 2019.
\newblock \doi{10.1080/03081079.2019.1585432}.

\bibitem[Boz{\'o}ki et~al.(2010)Boz{\'o}ki, F{\"u}l{\"o}p, and
  R{\'o}nyai]{Bozoki2010ooco}
S.~Boz{\'o}ki, J.~F{\"u}l{\"o}p, and L.~R{\'o}nyai.
\newblock On optimal completion of incomplete pairwise comparison matrices.
\newblock \emph{Mathematical and Computer Modelling}, 52\penalty0
  (1--2):\penalty0 318 -- 333, 2010.
\newblock ISSN 0895-7177.
\newblock \doi{http://dx.doi.org/10.1016/j.mcm.2010.02.047}.

\bibitem[Brunelli and Mezei(2017)]{Brunelli2017aiia}
M.~Brunelli and J.~Mezei.
\newblock {An inquiry into approximate operations on fuzzy numbers}.
\newblock \emph{International Journal of Approximate Reasoning}, 81:\penalty0
  147--159, 2017.
\newblock ISSN 0888613X.
\newblock \doi{10.1016/j.ijar.2016.11.011}.
\newblock URL \url{http://dx.doi.org/10.1016/j.ijar.2016.11.011}.

\bibitem[Brunelli and Rezaei(2019)]{Brunelli2019ambw}
M.~Brunelli and J.~Rezaei.
\newblock {A multiplicative best--worst method for multi-criteria decision
  making}.
\newblock \emph{Operations Research Letters}, 47\penalty0 (1):\penalty0 12--15,
  2019.
\newblock ISSN 01676377.
\newblock \doi{10.1016/j.orl.2018.11.008}.
\newblock URL \url{https://doi.org/10.1016/j.orl.2018.11.008}.

\bibitem[Cavallo and Brunelli(2018)]{Cavallo2018aguf}
B.~Cavallo and M.~Brunelli.
\newblock A general unified framework for interval pairwise comparison
  matrices.
\newblock \emph{International Journal of Approximate Reasoning}, 93:\penalty0
  178--198, 2018.
\newblock \doi{10.1016/j.ijar.2017.11.002}.

\bibitem[Cavallo and D'Apuzzo(2010)]{Cavallo2010cocp}
B.~Cavallo and L.~D'Apuzzo.
\newblock Characterizations of consistent pairwise comparison matrices over
  abelian linearly ordered groups.
\newblock \emph{Int. J. Intell. Syst.}, 25\penalty0 (10):\penalty0 1035--1059,
  2010.
\newblock \doi{10.1002/int.20438}.
\newblock URL \url{https://doi.org/10.1002/int.20438}.

\bibitem[Cavallo and D'Apuzzo(2012)]{Cavallo2012dwfa}
B.~Cavallo and L.~D'Apuzzo.
\newblock {Deriving weights from a pairwise comparison matrix over an
  alo-group}.
\newblock \emph{Soft Computing}, 16\penalty0 (2):\penalty0 353--366, 2012.
\newblock ISSN 14327643.
\newblock \doi{10.1007/s00500-011-0746-8}.

\bibitem[Cayley(1889)]{Cayley1889atot}
A.~Cayley.
\newblock A theorem on trees.
\newblock \emph{Quart. J. Pure Appl. Math.}, 23:\penalty0 376--378, 1889.

\bibitem[Chakraborty et~al.(2019)Chakraborty, Chowdhury, Chakraborty, Mehera,
  and Pal]{Chakraborty2019afga}
M.~Chakraborty, S.~Chowdhury, J.~Chakraborty, R.~Mehera, and R.~K. Pal.
\newblock {Algorithms for generating all possible spanning trees of a simple
  undirected connected graph: an extensive review}.
\newblock \emph{Complex and Intelligent Systems}, 5\penalty0 (3):\penalty0
  265--281, 2019.
\newblock ISSN 2199-4536.
\newblock \doi{10.1007/s40747-018-0079-7}.
\newblock URL \url{https://doi.org/10.1007/s40747-018-0079-7}.

\bibitem[Chang and Lee(1995)]{Chang1995teon}
P.~T. Chang and E.~S. Lee.
\newblock {The estimation of normalized fuzzy weights}.
\newblock \emph{Computers and Mathematics with Applications}, 29\penalty0
  (5):\penalty0 21--42, 1995.
\newblock ISSN 08981221.
\newblock \doi{10.1016/0898-1221(94)00246-H}.

\bibitem[Crawford and Williams(1985)]{Crawford1985anot}
R.~Crawford and C.~Williams.
\newblock A note on the analysis of subjective judgement matrices.
\newblock \emph{Journal of Mathematical Psychology}, 29:\penalty0 387 -- 405,
  1985.

\bibitem[Dede et~al.(2021)Dede, Kamalakis, and Anagnostopoulos]{Dede2021afoi}
G.~Dede, T.~Kamalakis, and D.~Anagnostopoulos.
\newblock {A framework of incorporating confidence levels to deal with
  uncertainty in pairwise comparisons}.
\newblock \emph{Central European Journal of Operations Research}, 2021.
\newblock ISSN 16139178.
\newblock \doi{10.1007/s10100-020-00735-0}.
\newblock URL \url{https://doi.org/10.1007/s10100-020-00735-0}.

\bibitem[Deng(1999)]{Deng1999mawf}
H.~Deng.
\newblock {Multicriteria analysis with fuzzy pairwise comparison}.
\newblock \emph{International Journal of Approximate Reasoning}, 21\penalty0
  (3):\penalty0 215--231, August 1999.

\bibitem[Diestel(2005)]{Diestel2005gt}
Reinhard Diestel.
\newblock \emph{{Graph theory}}.
\newblock Springer Verlag, 2005.

\bibitem[Dubois and Prade(1980)]{Dubois1980fsas}
D.~Dubois and H.~Prade.
\newblock \emph{Fuzzy sets and systems. Theory and applications}, volume 144 of
  \emph{Mathematics in Science and Engineering}.
\newblock Academic Press, 1980.

\bibitem[Durbach and Stewart(2012)]{Durbach2012muim}
I.~N. Durbach and T.~J. Stewart.
\newblock {Modeling uncertainty in multi-criteria decision analysis}.
\newblock \emph{European Journal of Operational Research}, 223\penalty0
  (1):\penalty0 1--14, 2012.
\newblock ISSN 03772217.
\newblock \doi{10.1016/j.ejor.2012.04.038}.
\newblock URL \url{http://dx.doi.org/10.1016/j.ejor.2012.04.038}.

\bibitem[Grabisch et~al.(2009)Grabisch, Marichal, Mesiar, and
  Pap]{Grabisch2009eoma}
M.~Grabisch, J.~L. Marichal, R.~Mesiar, and E.~Pap.
\newblock \emph{Aggregation functions}, volume 127 of \emph{Encyclopedia of
  mathematics and its applications}.
\newblock Cambridge University Press, Cambridge, England, July 2009.

\bibitem[Greco et~al.(2016)Greco, Ehrgott, and Figueira]{Greco2016mcda}
S.~Greco, M.~Ehrgott, and J.~Figueira, editors.
\newblock \emph{Multiple Criteria Decision Analysis: State of the Art Surveys}.
\newblock Springer, 2016.

\bibitem[Gross et~al.(2018)Gross, Yellen, and Anderson]{Gross2018gtaI}
J.~L. Gross, J.~Yellen, and M.~Anderson.
\newblock \emph{{Graph Theory and Its Applications}}.
\newblock Chapman and Hall/{CRC}, November 2018.
\newblock \doi{10.1201/9780429425134}.
\newblock URL \url{https://doi.org/10.1201/9780429425134}.

\bibitem[Harker(1987)]{Harker1987amoq}
P.~T. Harker.
\newblock Alternative modes of questioning in the analytic hierarchy process.
\newblock \emph{Mathematical Modelling}, 9\penalty0 (3):\penalty0 353 -- 360,
  1987.
\newblock ISSN 0270-0255.
\newblock \doi{https://doi.org/10.1016/0270-0255(87)90492-1}.

\bibitem[Janicki and Lenar{\v c}i{\v c}(2016)]{Janicki2016oawr}
R.~Janicki and A.~Lenar{\v c}i{\v c}.
\newblock Optimal approximations with rough sets and similarities in measure
  spaces.
\newblock \emph{International Journal of Approximate Reasoning}, 71:\penalty0
  1--14, 2016.
\newblock ISSN 0888-613X.
\newblock \doi{https://doi.org/10.1016/j.ijar.2015.12.014}.
\newblock URL
  \url{https://www.sciencedirect.com/science/article/pii/S0888613X15002029}.

\bibitem[Kaufmann and Gupta(1988)]{Kaufmann1988fmmi}
A.~Kaufmann and M.~M. Gupta.
\newblock \emph{Fuzzy Mathematical Models in Engineering and Management
  Science}.
\newblock Elsevier Science Publishers, North-Holland, Amsterdam, N.Y., 1988.

\bibitem[Kelmans and Chelnokov(1974)]{Kelmans1974acpo}
A.~K. Kelmans and V.~M. Chelnokov.
\newblock A certain polynomial of a graph and graphs with an extremal number of
  trees.
\newblock \emph{Journal of Combinatorial Theory (B)}, 16:\penalty0 197--214,
  1974.

\bibitem[Kirchhoff(1847)]{Kirchhoff1847udad}
G.~Kirchhoff.
\newblock Ueber der aulosung der gleichungen, auf welche man bei der
  untersuchung der linearen verteilung galvanischer strome gefuhrt wird.
\newblock \emph{Annalen der Physik}, 148\penalty0 (12):\penalty0 497--508,
  1847.

\bibitem[Ku{\l}akowski(2018)]{Kulakowski2018iito}
K.~Ku{\l}akowski.
\newblock Inconsistency in the ordinal pairwise comparisons method with and
  without ties.
\newblock \emph{European Journal of Operational Research}, 270\penalty0
  (1):\penalty0 314 -- 327, 2018.
\newblock ISSN 0377-2217.
\newblock \doi{https://doi.org/10.1016/j.ejor.2018.03.024}.

\bibitem[Ku{\l}akowski(2020{\natexlab{a}})]{Kulakowski2020utahp}
K.~Ku{\l}akowski.
\newblock \emph{Understanding the Analytic Hierarchy Process}.
\newblock Chapman and Hall / CRC Press, 2020{\natexlab{a}}.
\newblock ISBN 9781315392226.
\newblock \doi{10.1201/b21817}.

\bibitem[Ku{\l}akowski(2020{\natexlab{b}})]{kulakowski2020otgm}
K.~Ku{\l}akowski.
\newblock On the geometric mean method for incomplete pairwise comparisons.
\newblock \emph{Mathematics}, 8\penalty0 (11):\penalty0 1--12, October
  2020{\natexlab{b}}.

\bibitem[Ku{\l}akowski et~al.(2019)Ku{\l}akowski, Mazurek, Ram{\'\i}k, and
  Soltys]{Kulakowski2019witc}
K.~Ku{\l}akowski, J.~Mazurek, J.~Ram{\'\i}k, and M.~Soltys.
\newblock {When is the condition of order preservation met?}
\newblock \emph{European Journal of Operational Research}, 277\penalty0
  (1):\penalty0 248--254, August 2019.

\bibitem[Kurniawati and Yuliando(2015)]{Kurniawati2015pios}
D.~Kurniawati and H.~Yuliando.
\newblock {Productivity Improvement of Small Scale Medium Enterprises (SMEs) on
  Food Products: Case at Yogyakarta Province, Indonesia}.
\newblock \emph{Agriculture and Agricultural Science Procedia}, 3:\penalty0
  189--194, 2015.
\newblock ISSN 22107843.
\newblock \doi{10.1016/j.aaspro.2015.01.037}.
\newblock URL \url{http://dx.doi.org/10.1016/j.aaspro.2015.01.037}.

\bibitem[Lewin et~al.(2015)Lewin, Glenton, Munthe-Kaas, Carlsen, Colvin,
  Gulmezoglu, Noyes, Booth, Garside, and Rashidian]{Lewin2015uqei}
S.~Lewin, C.~Glenton, H.~Munthe-Kaas, B.~Carlsen, C.~J. Colvin, M.~Gulmezoglu,
  J.~Noyes, A.~Booth, R.~Garside, and A.~Rashidian.
\newblock {Using Qualitative Evidence in Decision Making for Health and Social
  Interventions: An Approach to Assess Confidence in Findings from Qualitative
  Evidence Syntheses (GRADE-CERQual)}.
\newblock \emph{PLoS Medicine}, 12\penalty0 (10):\penalty0 1--18, 2015.
\newblock ISSN 15491676.
\newblock \doi{10.1371/journal.pmed.1001895}.

\bibitem[Li and Ma(2020)]{Li2020crid}
H.~H. Li and W.~J. Ma.
\newblock {Confidence reports in decision-making with multiple alternatives
  violate the Bayesian confidence hypothesis}.
\newblock \emph{Nature Communications}, 11\penalty0 (1), 2020.
\newblock ISSN 20411723.
\newblock \doi{10.1038/s41467-020-15581-6}.
\newblock URL \url{http://dx.doi.org/10.1038/s41467-020-15581-6}.

\bibitem[Lundy et~al.(2016)Lundy, Siraj, and Greco]{Lundy2016tmeo}
M.~Lundy, S.~Siraj, and S.~Greco.
\newblock {The mathematical equivalence of the ``spanning tree'' and row
  geometric mean preference vectors and its implications for preference
  analysis}.
\newblock \emph{European Journal of Operational Research}, pages 1--12,
  September 2016.

\bibitem[Mardani et~al.(2015)Mardani, Jusoh, Nor, Khalifah, Zakwan, and
  Valipour]{Mardani2015mcdm}
A.~Mardani, A.~Jusoh, K.~M.D. Nor, Z.~Khalifah, N.~Zakwan, and A.~Valipour.
\newblock {Multiple criteria decision-making techniques and their applications
  - A review of the literature from 2000 to 2014}.
\newblock \emph{Economic Research-Ekonomska Istrazivanja}, 28\penalty0
  (1):\penalty0 516--571, 2015.
\newblock ISSN 1331677X.
\newblock \doi{10.1080/1331677X.2015.1075139}.
\newblock URL \url{http://dx.doi.org/10.1080/1331677X.2015.1075139}.

\bibitem[Mazurek and Ku{\l}akowski(2020)]{Mazurek2020sotc}
J.~Mazurek and K~Ku{\l}akowski.
\newblock {Satisfaction of the condition of order preservation: A simulation
  study.}
\newblock \emph{Operations Research and Decisions}, 30\penalty0 (2):\penalty0
  77--89, 2020.
\newblock \doi{10.37190/ord200205}.
\newblock URL
  \url{https://journals.indexcopernicus.com/search/article?articleId=2669489}.

\bibitem[Mazurek and Ram{\'\i}k(2019)]{Mazurek2019snpo}
J.~Mazurek and J.~Ram{\'\i}k.
\newblock {Some new properties of inconsistent pairwise comparisons matrices}.
\newblock \emph{International Journal of Approximate Reasoning}, 113:\penalty0
  119--132, October 2019.

\bibitem[Peterson and Pitz(1988)]{Peterson1988cuat}
D.~K. Peterson and G.~F. Pitz.
\newblock Confidence, uncertainty, and the use of information.
\newblock \emph{Journal of Experimental Psychology: Learning, Memory, and
  Cognition}, 14:\penalty0 85--92, 1988.

\bibitem[Quarteroni et~al.(2000)Quarteroni, Sacco, and
  Saleri]{Quarteroni2000nm}
A.~Quarteroni, R.~Sacco, and F.~Saleri.
\newblock \emph{Numerical mathematics}.
\newblock Springer Verlag, 2000.

\bibitem[Ramik(2020)]{Ramik2020pcmt}
J.~Ramik.
\newblock \emph{Pairwise Comparisons Method: Theory and Applications in
  Decision Making}.
\newblock Springer, 2020.

\bibitem[Saaty(1977)]{Saaty1977asmf}
T.~L. Saaty.
\newblock A scaling method for priorities in hierarchical structures.
\newblock \emph{Journal of Mathematical Psychology}, 15\penalty0 (3):\penalty0
  234 -- 281, 1977.
\newblock ISSN 0022-2496.
\newblock \doi{10.1016/0022-2496(77)90033-5}.

\bibitem[Saaty(2004)]{Saaty2004dmta}
T.~L. Saaty.
\newblock {Decision making --- the Analytic Hierarchy and Network Processes
  (AHP/ANP)}.
\newblock \emph{Journal of Systems Science and Systems Engineering},
  13\penalty0 (1):\penalty0 1--35, 2004.
\newblock ISSN 1004-3756.
\newblock \doi{10.1007/s11518-006-0151-5}.
\newblock URL \url{http://dx.doi.org/10.1007/s11518-006-0151-5}.

\bibitem[Siraj et~al.(2012)Siraj, Mikhailov, and Keane]{Siraj2012east}
S.~Siraj, L.~Mikhailov, and J.~A. Keane.
\newblock {Enumerating all spanning trees for pairwise comparisons}.
\newblock \emph{Computers and Operations Research}, 39\penalty0 (2):\penalty0
  191--199, February 2012.

\bibitem[Tone(1993)]{Tone1993llsm}
K.~Tone.
\newblock {Logarithmic Least Squares Method for Incomplete Pairwise Comparisons
  in the Analytic Hierarchy Process}.
\newblock Technical Report 94-B-2, Saitama University, Institute for Policy
  Science Research, Urawa, Saitama, 338, Japan, December 1993.

\bibitem[Vaidya and Kumar(2006)]{Vaidya2006ahpa}
O.~S. Vaidya and S.~Kumar.
\newblock {Analytic hierarchy process: An overview of applications}.
\newblock \emph{European Journal of Operational Research}, 169\penalty0
  (1):\penalty0 1--29, February 2006.

\bibitem[Weiss and Rao(1987)]{Weiss1987adfl}
E.~N. Weiss and V.~R. Rao.
\newblock {AHP DesignIssues for Large-scale Systems}.
\newblock \emph{Decision Sciences}, 18\penalty0 (1):\penalty0 43--61, 1987.
\newblock ISSN 15405915.
\newblock \doi{10.1111/j.1540-5915.1987.tb01502.x}.

\bibitem[Wesson and Pulford(2009)]{Wesson2009veoc}
C.~J. Wesson and B.~D. Pulford.
\newblock {Verbal expressions of confidence and doubt}.
\newblock \emph{Psychological Reports}, 105\penalty0 (1):\penalty0 151--160,
  2009.
\newblock ISSN 00332941.
\newblock \doi{10.2466/PR0.105.1.151-160}.

\bibitem[William-West and Ciucci(2021)]{WilliamWest2021dtfw}
T.~O. William-West and D.~Ciucci.
\newblock {Decision-theoretic five-way approximation of fuzzy sets}.
\newblock \emph{Information Sciences}, 572:\penalty0 200--222, 2021.
\newblock ISSN 00200255.
\newblock \doi{10.1016/j.ins.2021.04.105}.
\newblock URL \url{https://doi.org/10.1016/j.ins.2021.04.105}.

\bibitem[Yang et~al.(2013)Yang, Yan, and Zeng]{Yang2013hthc}
X.~Yang, L.~Yan, and L.~Zeng.
\newblock {How to handle uncertainties in AHP: The Cloud Delphi hierarchical
  analysis}.
\newblock \emph{Information Sciences}, 222:\penalty0 384--404, 2013.
\newblock ISSN 00200255.
\newblock \doi{10.1016/j.ins.2012.08.019}.
\newblock URL \url{http://dx.doi.org/10.1016/j.ins.2012.08.019}.

\bibitem[Zadeh(1965)]{Zadeh1965fs}
L.~A. Zadeh.
\newblock {Fuzzy Sets}.
\newblock \emph{Information and Control}, 8:\penalty0 338--353, 1965.
\newblock \doi{10.1016/S0019-9958(65)90241-X}.

\bibitem[Zadeh(2011)]{Zadeh2011anoz}
L.~A. Zadeh.
\newblock {A Note on Z-numbers}.
\newblock \emph{Information Sciences}, 181\penalty0 (14):\penalty0 2923--2932,
  2011.
\newblock ISSN 00200255.
\newblock \doi{10.1016/j.ins.2011.02.022}.
\newblock URL \url{http://dx.doi.org/10.1016/j.ins.2011.02.022}.

\bibitem[Zavadskas et~al.(2014)Zavadskas, Turskis, and
  Kildiene]{Zavadskas2014soas}
E.~K. Zavadskas, Z.~Turskis, and S.~Kildiene.
\newblock {State of art surveys of overviews on MCDM/MADM methods}.
\newblock \emph{Technological and Economic Development of Economy}, 20\penalty0
  (1):\penalty0 165--179, 2014.
\newblock ISSN 20294921.
\newblock \doi{10.3846/20294913.2014.892037}.

\end{thebibliography}

\end{document}